\documentclass[10pt,a4paper]{amsart} 
\usepackage{graphicx,latexsym,amssymb,color}
\theoremstyle{plain} 
\newtheorem{thm}{Theorem}[section] 
\newtheorem{lem}[thm]{Lemma} 
\newtheorem{prop}[thm]{Proposition} 
\newtheorem{cor}[thm]{Corollary} 
\theoremstyle{definition}
 
\theoremstyle{remark} 
\newtheorem{rem}[thm]{Remark}

\numberwithin{equation}{section}
\numberwithin{figure}{section}
\newcommand{\bd}{\begin{description}}   
\newcommand{\ed}{\end{description}} 
\newcommand{\bF}{\begin{figure}[!h]}      \newcommand{\eF}{\end{figure}} 
\newcommand{\ba}{\begin{array}}      \newcommand{\ea}{\end{array}} 
\newcommand{\bc}{\begin{center}}     \newcommand{\ec}{\end{center}} 
\newcommand{\be}{\begin{enumerate}}  \newcommand{\ee}{\end{enumerate}} 
\newcommand{\beq}{\begin{eqnarray}}  \newcommand{\eeq}{\end{eqnarray}} 
\newcommand{\beQ}{\begin{eqnarray*}} \newcommand{\eeQ}{\end{eqnarray*}} 
\newcommand{\bi}{\begin{itemize}}    \newcommand{\ei}{\end{itemize}}

\newcommand{\ov}{\overline}

\begin{document} 
\title[Borromean surgery formula for the Casson invariant]{Borromean surgery formula for the Casson invariant} 
\author[J.B. Meilhan]{Jean-Baptiste Meilhan} 
\address{CTQM -- Department of Mathematical Sciences\\
         University of Aarhus \\
         Ny Munkegade, bldg. 1530, \\
         8000 Aarhus C, Denmark}
	 \email{meilhan@imf.au.dk}
\begin{abstract} 
It is known that every oriented integral homology $3$-sphere can be obtained from $S^3$ by  a finite sequence of Borromean surgeries.  
We give an explicit formula for the variation of the Casson invariant 
under such a surgery move.  The formula involves simple classical invariants, namely the framing, linking number and Milnor's triple linking number.  
A more general statement, for $n$ independent Borromean surgeries, is also provided.  
\end{abstract} 
\maketitle 
\section{Introduction} \label{intro}
A \emph{Borromean surgery link} in a closed oriented $3$-manifold $M$ is an oriented framed link obtained by embedding in $M$ a copy of the standard genus $3$ handlebody containing the $6$-component $0$-framed oriented link depicted below.  (Here and throughout this paper, blackboard framing convention is used.)  Surgery along such a link is called a \emph{Borromean surgery move}.  
\bc
\includegraphics{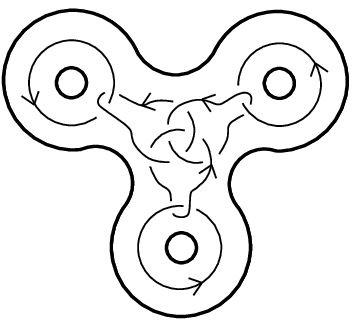}
\ec
The notion of Borromean surgery was first introduced by Matveev \cite{matv} in the eighties, and was more recently used to define the Goussarov-Habiro finite type invariant theory for closed oriented $3$-manifolds \cite{GGP,habi}.  

Matveev showed that two closed oriented 3-manifolds are \emph{Borromean equivalent}, i.e. are related by a sequence of 
Borromean surgery moves, if and only if they have the same first homology groups and isomorphic linking forms.  
In particular, any oriented integral homology $3$-sphere is obtained from $S^3$ by a finite sequence of Borromean surgeries.  
It is thus a natural problem to study the behavior of classical invariants of integral homology spheres, such as the Casson invariant, 
under this kind of surgery move.  
Our main result is Theorem \ref{cy1}, which expresses the variation of the Casson invariant under a Borromean surgery move 
in terms of simple classical invariants of the Borromean surgery link, namely framings, linking numbers and 
Milnor's triple linking number -- see \S \ref{statements} for an explicit statement.  
A more general formula, for $n$ independent Borromean surgeries on an integral homology sphere, is also provided in Theorem \ref{cy}, 
which involves an additional cubic expression in some linking numbers of the surgery link.    

The paper is organized as follows.  In \S \ref{S2}, we recall the definition and basic properties of the Casson invariant and state the above-mentioned surgery formulas.  
In \S \ref{S3} we prove a couple of preliminary results on Borromean surgery in $S^3$ using a crossing change formula for the Casson invariant due to Johannes \cite{johannes}.  
In \S \ref{proofL} we prove the formula of Theorem \ref{cy1} for a single Borromean surgery using the results of \S \ref{S3} and Lescop's sum formula for the Casson invariant \cite{L}.  
In \S \ref{proofzhs} we use the theory of finite type invariants to prove the formula of Theorem \ref{cy} for multiple Borromean surgeries.  
In \S \ref{rencom} we compare Theorem \ref{cy1} with a similar formula derived from Lescop's global surgery formula, and we apply it to a realization problem for homology spheres of Mazur type.  

\noindent \subsection*{Acknowledgments}
This paper benefited greatly from discussions with Christine Lescop at various stages of its writing.  
The author also thanks Kazuo Habiro, Gw\'e\-na\"el Massuyeau and Alex James Bene for helpful comments, and Christian Blanchet for suggesting working on this subject.  
\section{Variation of the Casson invariant under Borromean surgery} \label{S2}
\subsection{The Casson invariant of integral homology spheres} \label{seccasson}
In this paper, an \emph{integral homology sphere} will always mean an oriented integral homology $3$-sphere.
We denote by $\mathbf{Z}HS$ the set of integral homology spheres, considered up to orientation-preserving diffeomorphisms.  
\begin{thm}[Casson] \label{casson}
There exists a unique function 
\[ \lambda : \mathbf{Z}HS \longrightarrow \mathbf{Z} \]
such that, for every $M,N\in \mathbf{Z}HS$ and for every knot $K$ in $M$ :
\begin{enumerate}
\item[(1)] $\lambda(S^3)=0$.  
\item[(2)]  For $n\in \mathbf{Z}$, let $M_{K_n}$ be the result of $\frac{1}{n}$-Dehn surgery on $M$ along $K$.  Then:
\[ \lambda(M_{K_{n+1}})-\lambda(M_{K_n})=\frac{1}{2}\Delta''_K(1), \]
where $\Delta_K(t)$ denotes the Alexander polynomial of $K$.  
\end{enumerate}
Furthermore, 
\begin{enumerate}
\item[(3)]  $\lambda(M\sharp N)=\lambda(M)+\lambda(N)$, where $\sharp$ denotes the connected sum.
\item[(4)]  The Rochlin invariant $\mu$ is the mod 2 reduction of $\lambda$.  More precisely, 
$\mu(M)\equiv \lambda(M)\textrm{ (mod $2$)}$.  
\end{enumerate}
\end{thm}

\noindent This unique function is called the \emph{Casson invariant} of integral homology spheres.  
This invariant counts, in some sense, the conjugacy classes of irreducible representations of $\pi_1(M)$ in $SU(2)$.  

In \cite{W}, Walker extended the Casson invariant to a $\mathbf{Q}$-valued invariant of rational homology spheres, called the Casson-Walker invariant.    
In this paper, we will use the following normalization of the Casson-Walker invariant
\[  \lambda = \frac{\lambda_W}{2},  \] 
where $\lambda_W$ denotes the normalization adopted by Walker in \cite{W}.  
(Note that our notations are consistent, in the sense that if we restrict to integral homology spheres then $ \frac{\lambda_W}{2}$ coincides with Casson's original invariant.)  

The Casson-Walker invariant was extended by Lescop to all closed oriented $3$-manifolds in \cite{LG}.  
\subsection{Conventions} \label{not}  
Let $L$ be a Borromean surgery link in a rational homology sphere $M$.  

Up to isotopy one can always assume that there is a $3$-ball $B$ in $M$ which intersects $L$ as depicted 
below.  The boundary of $B$ intersects three components of $L$, called the \emph{leaves} of $L$.  
\begin{figure}[!h]
\input{conv_new.pstex_t}
\caption{} \label{conv}
\end{figure}
\noindent 
We fix an order on the set of the leaves of $L$, and we denote them by $F_1$, $F_2$ and $F_3$ according to this order.  
Denote by $f_i$ the framing of $F_i$, and by $l_{ij}$ the linking number $lk(F_i,F_j)$, $1\le i\ne j\le 3$.   

\begin{figure}[!h]
\input{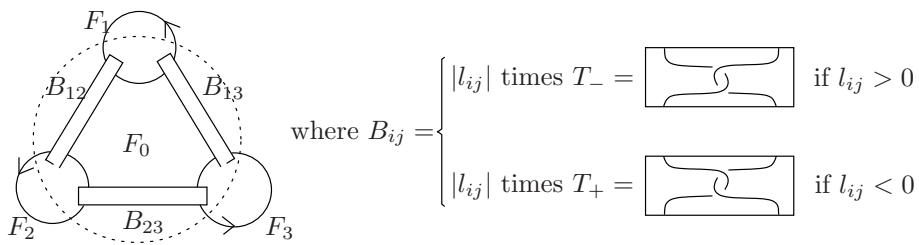}
\caption{The algebraically split link $F_0$. } \label{F0}
\end{figure}

Observe that there is a canonical $3$-component algebraically split link $F_0$ associated to the Borromean surgery link $L$ as follows.  
For each pair $1\le i<j\le 3$, consider the $(2,2)$-tangle $B_{ij}$ obtained by stacking  $|l_{ij}|$ copies of the clasp-shaped $(2,2)$-tangle 
$T_{\pm}$ depicted on the right-hand side of Figure \ref{F0} (depending on the sign of $l_{ij}$).  
If $l_{ij}=0$, $B_{ij}$ is the trivial $(2,2)$-tangle.  
Starting with the link $F_1\cup F_2\cup F_3$, insert for each $1\le i<j\le 3$ the tangle $B_{ij}$ into the leaves $F_i$ and $F_j$ 
so that their linking number is reduced to zero, see Figure \ref{F0}.  
The resulting $3$-component link is the algebraically split link $F_0$ associated to $L$.  

Because $F_0$ has linking numbers zero, it has a well-defined \emph{Milnor's triple linking number}  $\ov{\mu}_{ijk}(F_0)$ in $M$ (a short definition is given in \S \ref{123} at the end of this section).  
In the following, we denote by $\mu_{123}(L;M)$ the triple linking number of the algebraically split link $F_0$ associated to $L$ in $M$.  We will usually simply write $\mu_{123}(L)$ when $M$ is explicit from the context.  
\subsection{Statement of the main results} \label{statements}  
We can now state our two main results.

We first express the variation of the Casson-Walker invariant under surgery along a Borromean surgery link $L$ in a rational homology sphere $M$, under certain restrictions.  
\begin{thm} \label{cy1}
Let $L$ be a Borromean surgery link in a rational homology sphere $M$, such that the leaves of $L$ have integral linking numbers and have self-linking numbers zero in $\mathbf{Q} / \mathbf{Z}$.  
Denote by $M_L$ the result of surgery on $M$ along $L$.  Then the difference $\lambda(M_L)-\lambda(M)$ is given by the formula
\[ - f_1.f_2.f_3 - 2.l_{12}.l_{13}.l_{23} - 2.\mu_{123}(L) 
+ \sum_{\circlearrowleft_{1,2,3}} l_{23}.(l_{23}+1).f_1, \]
\noindent where the sum is over all cyclic permutations of the indices $(1,2,3)$.  
\end{thm}

Let us now consider the case of a disjoint union $L_1\cup...\cup L_n$ of $n$ Borromean surgery links in an integral homology sphere $M$.  
We use the conventions of \S \ref{not} for the orientation and ordering of the leaves of each link: 
the leaves of $L_k$ are denoted by $F^k_1$, $F^k_2$ and $F^k_3$ ($1\le k\le n$).  
Denote by $l^{kl}_{ij}$ the linking number $lk(F^k_i,F^l_j)$ ; $1\le i\ne j\le 3$, $1\le k\ne l\le n$.   
\begin{thm} \label{cy}
The difference $\lambda(M_{L_1\cup...\cup L_n})-\lambda(M)$ is given by the formula
\[ \sum_{1\le i\le n} \left(\lambda(M_{L_i}) -\lambda(M)\right) - 
2\sum_{1\le k<l \le n} \sum_{\sigma\in  S_3} s(\sigma) l^{kl}_{1\sigma(1)}l^{kl}_{2\sigma(2)}l^{kl}_{3\sigma(3)}, \]
\noindent where, for an element $\sigma$ in the symmetric group $S_3$, $s(\sigma)\in\{-1;+1\}$ denotes its signature.  
\end{thm}

By reducing the formula in Theorem \ref{cy} modulo 2 and using Theorem \ref{casson}(4), we obtain a similar formula for the variation of Rochlin's 
$\mu$-invariant of an integral homology sphere $M$ under surgery along a disjoint union $L_1\cup ...\cup L_n$ of Borromean surgery links.  
\begin{cor}
\[ \mu(M_{L_1\cup ...\cup L_n}) - \mu(M) = \sum_{1\le i\le n} f^i_1.f^i_2.f^i_3 \textrm{ (mod $2$)}.  \]
\end{cor}
This formula was already (essentially) established independently by Massuyeau in \cite{massuyeau}, and Auclair and Lescop in \cite{AL}.  
\subsection{Milnor's triple linking number}\label{123}
Let $K=K_1\cup K_2\cup K_3$ be an algebraically split, oriented, ordered link in a rational homology sphere $M$.  
Pick a triple $(S_1,S_2,S_3)$ of transverse Seifert surfaces such that each $S_i$ is bounded by $K_i$ ($i=1,2,3$) and does not intersect the two other components of $K$ (this is possible since all linking numbers are zero).  Then \emph{Milnor's triple linking number} $\ov{\mu}_{123}(K;M)$ of $K$ in $M$ is defined as the negative of the algebraic intersection $\langle S_1,S_2,S_3\rangle_M$, defined as follows.  

Each surface $S_i$ is equipped with a positive normal vector field, induced by its orientation and the orientation of $M$, and thus each triple intersection point $x\in S_1\cap S_2\cap S_3$ can be associated a sign, depending on whether the ordered basis of normal vectors to $S_1$, $S_2$ and $S_3$ at $x$ agrees or not with the orientation of $M$. The algebraic intersection $\langle S_1,S_2,S_3\rangle_M$ is the sum of these signs over all $x\in S_1\cap S_2\cap S_3$.  
\section{Preliminary results: Borromean surgery in $S^3$} \label{S3}
In this section we prove a couple of preliminary lemmas on Borromean surgery in $S^3$.  
\subsection{Main tools}
We recall here a self-crossing change formula for the Casson invariant due to Johannes, as well as a simple result derived from Kirby calculus.  
\subsubsection{Johannes' formula for the Casson invariant} \label{J}
In \cite{johannes}, Johannes expresses the difference of the Casson invariants of $3$-manifolds presented by 
links in $S^3$ which differ by a crossing change within a component.  
This type of move on surgery links was studied in \cite{matv} under the name of \emph{Whitehead surgery}.   

\begin{figure}[!h] 
\includegraphics{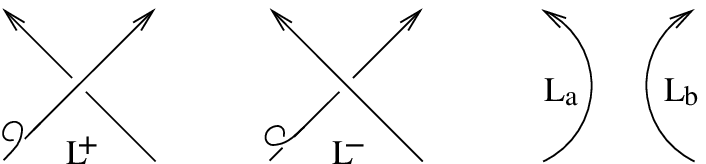}
\caption{} \label{cross}
\end{figure}

For the purpose of the present paper, it is enough to consider the case of $2$-component surgery links.  
Let $L^+=L^+_1\cup L_2$ and $L^-=L^-_1\cup L_2$ be two framed links in $S^3$, with same framing, 
which only differ by a crossing change on the first component.  
Let $L_a\cup L_b$ be the $2$-component link obtained from $L^{\pm}_1$ by smoothing this crossing (see Figure \ref{cross}).
Then Johannes' formula \cite{johannes} states that  
\begin{equation} \label{johannes}
 \lambda(S^3_{L^+})-\lambda(S^3_{L^-})=\dfrac{f_2.l_{ab} - l_{a2}.l_{b2}}{f_1.f_2 - l_{12}^2},
\end{equation}
where $f_1$ (resp.\ $f_2$) denotes the framing of $L^\pm_1$ (resp.\ $L_2$), 
and $l_{ab}$ (resp.\ $l_{a2}$, $l_{b2}$, $l_{12}$) denotes the linking number $lk(L_a,L_b)$ (resp.\ $lk(L_a,L_2)$, $lk(L_b,L_2)$, $lk(L^{\pm}_1,L_2)$).   
Observe that the denominator in (\ref{johannes}) is merely the determinant of the linking matrix of $L^{\pm}$.  

Note that the formulas given in \cite{johannes} are for the Casson-Walker invariant, using Walker's normalization (see \S \ref{seccasson}).  Note also that 
we fixed a sign error in \cite{johannes} so that (\ref{johannes}) agrees with Theorem \ref{casson}.\footnote{For example, 
consider the right-handed trefoil knot $T$ in $S^3$. It can be easily verified by Theorem \ref{casson}(2) that $\lambda(S^3_{T_{+1}})=1$, 
but \cite{johannes} would give $-1$.  }
\begin{rem}
Surgery along a $\pm 1$-framed knot and Whitehead surgery are two kinds of surgery moves relating all integral homology spheres, and   
Casson and Johannes respectively gave a formula for the variation of the Casson invariant under these moves.  
In this paper, we give such a formula for a third type of move which relates all integral homology spheres, namely 
Borromean surgery.  
\end{rem}
\subsubsection{Borromean surgery link and Kirby equivalence} \label{Ki}
In addition to Johannes' formula, we will need the following lemma, which follows directly from several applications of Kirby moves (see \cite{GGP,habi}).  
\begin{lem} \label{KIRBY}
Let $K$ be a Borromean surgery link.  
$K$ is Kirby-equivalent to the $2$-component link $\tilde{K}$ represented in Figure \ref{kirby}.  
\bF
\input{Kir.pstex_t}
\caption{Two Kirby-equivalent links.} \label{kirby}
\eF
\end{lem}
Observe that the linking number of $\tilde{K}$ is $1$, and that one of its component has always framing $0$.  The determinant of the linking matrix of $\tilde{K}$ is thus always $-1$.  
\subsection{Two lemmas in $S^3$}\label{lems}
We can now establish the two above-mentioned preliminary lemmas.  
In both statements, we consider a local operation on a Borromean surgery link in $S^3$, and give a formula for the variation of the Casson invariant of the surgered manifolds.  
%
%
\subsubsection{Changing the framing.}
We first consider the operation of adding a kink to a leaf of a Borromean surgery link.  
\begin{lem}\label{lemfr}
Let $K$ be a Borromean surgery link in $S^3$, and let $K_i$ be obtained by changing the framing of the $i^{th}$ leaf of $K$ by $+1$ as depicted below.  
\bc
\includegraphics{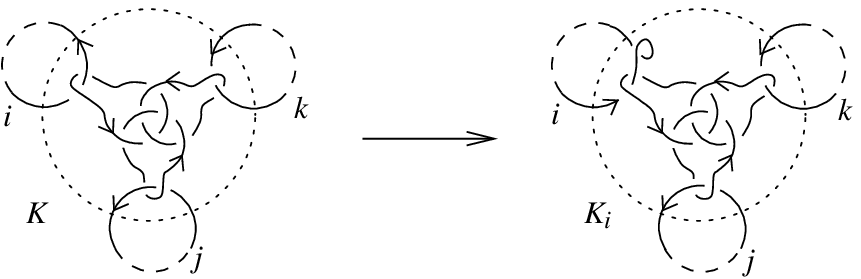}
\ec
We have: 
$$\lambda(S^3_K)= \lambda(S^3_{K_i}) - f_j.f_k + l_{jk}(1+l_{jk}).  $$
\end{lem}
\begin{proof}
By Lemma \ref{KIRBY}, we have $\lambda(S^3_K)=\lambda(S^3_{\tilde{K}})$ and $\lambda(S^3_{K_i})=\lambda(S^3_{\tilde{K}_i})$, 
where $\tilde{K}$ and $\tilde{K}_i$ are two $2$-component links as depicted in Figure \ref{linkfr}.  
\bF
\input{frlinknew.pstex_t}
\caption{} \label{linkfr}
\eF
Consider the $2$-component link $K_1\cup K_2$ isotopic to $\tilde{K}_i$, and the positive crossing $c$ in $K_1\cup K_2$ as in Figure \ref{linkfr}.  Changing this crossing (without changing the framing) gives the link $\tilde{K}$, and smoothing this crossing turns $K_1$ into a $2$-component link $L_a\cup L_b$ satisfying $lk(L_a,L_b)=-f_k$, $lk(L_a,K_2)=1+l_{jk}$ and $lk(L_b,K_2)=-l_{jk}$.  
The result then follows from Johannes' formula.  
\end{proof}
\subsubsection{Changing the linking numbers} 
We now study the effect, at the level of the Casson invariant, of adding a clasp between two leaves of a Borromean surgery link.  
\begin{lem}\label{lemlk}
Let $K$ be a Borromean surgery link in $S^3$, and let $K^-$ be obtained by adding a clasp between the leaves $i$ and $j$ of $K$ as shown below.  
\bc
\includegraphics{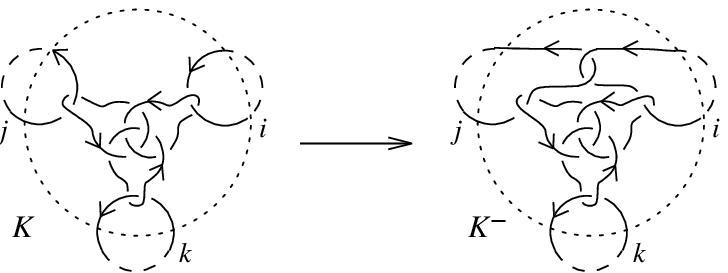}
\ec
We have: 
$$\lambda(S^3_K)= \lambda(S^3_{K^-})-2.l_{ik}.l_{jk} + 2.f_k.l_{ij}.$$  
\end{lem}
\begin{proof}
By Lemma \ref{KIRBY}, we have $\lambda(S^3_K)=\lambda(S^3_{\tilde{K}})$ and $\lambda(S^3_{K^-})=\lambda(S^3_{\tilde{K}^-})$, 
where $\tilde{K}$ and $\tilde{K}^-$ are two $2$-component links as depicted in Figure \ref{linkfr} and \ref{K} respectively.  Now observe that 
$\tilde{K}^-$ is isotopic to the link $K_1\cup K_2$ of Figure \ref{K}, which only differs from $\tilde{K}$ by four crossing changes within the component $K_1$.  
\bF
\input{linknew.pstex_t}
\caption{} \label{K}
\eF 
The result thus follows from $4$ applications of Johannes' formula, in a same way as in the previous proof.
\end{proof}
\section{Proof of Theorem \ref{cy1}} \label{proofL}
We can now prove Theorem \ref{cy1}.  Our main tool here is Lescop's sum formula for the Casson invariant \cite{L}, which we briefly review in the next subsection.   
\subsection{Lescop's sum formula} \label{sum}
Given a compact $3$-manifold $M$ with boundary, we denote by $\mathcal{L}_M$ the kernel of the map 
$H_1(\partial M;\mathbf{Q})\rightarrow H_1(M;\mathbf{Q})$ induced by the inclusion.  We call $\mathcal{L}_M$ the \emph{Lagrangian} of $M$.  

Let $\Sigma$ be a closed oriented surface of genus $g$, and let $A$, $B$, $A'$, $B'$ be four rational homology handlebodies such that 
$\partial A=\partial A'=-\partial B=-\partial B'=\Sigma$. 
Assume that in $H_1(\Sigma;\mathbf{Q})$ we have $\mathcal{L}_{A}=\mathcal{L}_{A'}$, $\mathcal{L}_{B}=\mathcal{L}_{B'}$, and $\mathcal{L}_A\cap \mathcal{L}_B =0$.  
Then the main result of \cite{L} states that $\lambda (A\cup B)-\lambda (A '\cup B)-\lambda (A\cup B ')+\lambda (A '\cup B ')$ is given by $$ -2\sum_{\{i,j,k\}\subset \{1,...,g\}}\mathcal{I}_{AA'}(\alpha_i\wedge \alpha_j\wedge \alpha_k) \mathcal{I}_{BB'}(\beta_i\wedge \beta_j\wedge \beta_k), $$ 
\noindent where $\{\alpha_i \}_i$, resp.\ $\{\beta_i \}_i$, is a basis for $\mathcal{L}_{A}$, resp.\ $\mathcal{L}_{B}$, such that the intersection form 
$\cdot$ on $\Sigma$ satisfies $\alpha_i\cdot \beta_j=\delta_{ij}$, and where $\mathcal{I}_{AA'}$ is the \emph{intersection form on $\Lambda^3 \mathcal{L}_{A}$} defined as follows.     
There is an isomorphism $H_2(A\cup_{\Sigma} -A';\mathbf{Q})\rightarrow \mathcal{L}_{A}$ coming from the Mayer-Vietoris sequence which, to the homology class of an oriented surface $S\subset A\cup_{\Sigma} -A'$, associates the class of the (oriented) boundary of $S\cap A$.  
Then $\mathcal{I}_{AA'}$ is the trilinear alternating form on $\mathcal{L}_{A}$ induced, \emph{via} this isomorphism, by the triple intersection form on $H_2(A\cup_{\Sigma} -A';\mathbf{Q})$.  

Note that the above sum formula is for the Casson-Walker invariant, with the normalization specified in \S \ref{seccasson}. 
\subsection{Proof of Theorem \ref{cy1}}
Let $L$ be a Borromean surgery link in a rational homology sphere $M$ such that the linking numbers of the leaves are integers $l_{ij}$, and such that the leaves have self-linking numbers zero in $\mathbf{Q} / \mathbf{Z}$, and thus have integral framings $f_i$.  
We aim to compute $\lambda(M_L)- \lambda(M)$.  

\begin{figure}[!h]
\input{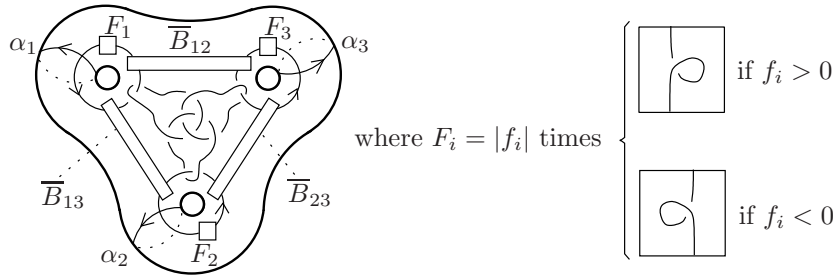}
\caption{The link $L_{st}$ in the handlebody $A$.  Here $\ov{B}_{ij}$ is $|l_{ij}|$ times the $(2,2)$-tangle $T_+$ (resp.\ $T_-$) of Figure \ref{F0} if $l_{ij}>0$ (resp.\ $l_{ij}<0$).  }
\label{figA}
\end{figure}

Consider the Borromean surgery link $L_{st}$ of Figure \ref{figA} in the genus three handlebody $A$, standardly embedded in $S^3$.  Denote by $B$ 
the exterior of $A$ in $S^3$, and set $A'=A_{L_{st}}$.  
The link $L$ can be regarded as obtained by an embedding of $A$ in $M$ such that the Lagrangians of $B$ and $B'=M\setminus \textrm{Int}(A)$ (\emph{i.e.} the curves of $\partial A$ that bound in $B$ and in $B'$) are the same.  
The main theorem of \cite{L} gives  
$$ \lambda(S^3) - \lambda(S^3_{L_{st}}) - \lambda(M) + \lambda(M_L)
 = -2 \mathcal{I}_{AA'}(\alpha_1\wedge \alpha_2\wedge \alpha_3) \mathcal{I}_{BB'}(\beta_1\wedge \beta_2\wedge \beta_3), $$
where the $\alpha_i$ are meridians of the three leaves of $L_{st}$ (see Figure \ref{figA}), and where the $\beta_i$ are longitudes of the three leaves of the link $F_0$ of Figure \ref{F0}.  

For $i=1,2,3$, let $\Sigma_{B}^i$ (resp.\ $\Sigma_{B'}^i$) be oriented surfaces in $B$ (resp.\ $B'$) bounded by $\beta_i$, such that all surfaces are transverse to each other and to $\partial B$ in $B\cup_{\partial B} -B'$.  Denote by $\Sigma_i$ the closed surface $\Sigma_{B}^i\cup_{\beta_i}\Sigma_{B'}^i$ in $B\cup_{\partial B} -B'$. Using the notations of \S \ref{123}, we have
$$ \mathcal{I}_{BB'}(\beta_1\wedge \beta_2\wedge \beta_3) 
 = \langle \Sigma_1,\Sigma_2,\Sigma_3\rangle_{B\cup_{\partial B} -B'}
 = \langle \Sigma_{B}^1,\Sigma_{B}^2,\Sigma_{B}^3\rangle_{B} - \langle \Sigma_{B'}^1,\Sigma_{B'}^2,\Sigma_{B'}^3\rangle_{B'},  $$ 
where $\langle \Sigma_{B}^1,\Sigma_{B}^2,\Sigma_{B}^3\rangle_{B}$ is zero, and where $- \langle \Sigma_{B'}^1,\Sigma_{B'}^2,\Sigma_{B'}^3\rangle_{B'}$ identifies directly with $\mu_{123}(L;M)$.  It also follows from the above and a direct computation that $\mathcal{I}_{AA'}(\alpha_1\wedge \alpha_2\wedge \alpha_3)=1$.  We thus obtain 
$$ \lambda(M_L) - \lambda(M) = -2\mu_{123}(L) +  \lambda(S^3_{L_{st}}). $$

Now, we can use the two lemmas of \S \ref{S3} to compute $\lambda(S^3_{L_{st}})$ as follows.  
By applying Lemma \ref{lemfr} $|f_1|+|f_2|+|f_3|$ times to the link $L_{st}$, we can reduce the framings of its three leaves to zero, and the resulting Borromean surgery link $L'$ satisfies 
\[ \lambda(S^3_{L_{st}})=\lambda(S^3_{L'}) - f_1.f_2.f_3 + \sum_{\circlearrowleft_{1,2,3}} f_1.l_{23}.(l_{23}+1). \]
Next, by suitably applying Lemma \ref{lemlk} $|l_{23}|+|l_{13}|+|l_{12}|$ times to the link $L'$, we obtain a Borromean surgery link $L_f$ in $S^3$ whose leaves are three $0$-framed pairwise unlinked unknots, and such that  
\[ \lambda(S^3_{L'})=\lambda(S^3_{L_f})-2.l_{12}.l_{13}.l_{23}. \]
Using Kirby calculus, one can easily check that $S^3_{L_f}\cong S^3$ (see for example \cite[\S 2.1]{GGP}).  
We thus have $\lambda(S^3_{L_f}) = 0$.  It follows that   
$$ \lambda(S^3_{L_{st}})=-f_1.f_2.f_3 - 2.l_{12}.l_{13}.l_{23} + \sum_{\circlearrowleft_{1,2,3}} l_{23}.(l_{23}+1).f_1, $$ 
which concludes the proof.  
\section{Proof of Theorem \ref{cy}} \label{proofzhs}
In order to prove Theorem \ref{cy}, we will use the theory of finite type invariants.  
\subsection{Finite type invariants of integral homology spheres} \label{Sfti}
The notion of finite type invariants of integral homology spheres is due to Ohtsuki, and is defined using algebraically split, unit-framed links \cite{O}.  For the purpose of the present paper, we use an alternative definition using the notion of Borromean surgery, due to Goussarov and Habiro \cite{G,GGP,habi}.

Let $\mathcal{S}$ denote the free $\mathbf{Z}$-module generated by elements of $\mathbf{Z}HS$.  
For $k\ge 0$, let $\mathcal{S}_k$ denote the $\mathbf{Z}$-submodule of $\mathcal{S}$ generated by elements 
  $$ [M;G]:=\sum_{G'\subseteq G} (-1)^{|G'|} M_{G'}, $$
where $M$ is an integral homology sphere, and where $G=G_1\cup...\cup G_k$ is a disjoint union of $k$ Borromean surgery links in $M$.  The sum runs over all the subsets $G'$ of $G$ (regarded as the set of the links $G_i$) and $|G'|$ denotes the cardinality of $G'$.  
Observe that we have the equality 
\begin{equation} \label{MG}
M_G=\sum_{G'\subseteq G } (-1)^{|G'|} [M;G'].  
\end{equation}

A \emph{finite type invariant of degree $k$} is a map $f: \mathbf{Z}HS\longrightarrow A$, where $A$ is an abelian group, whose natural extension to $\mathcal{S}$ vanishes on $\mathcal{S}_{k+1}$.  
\begin{prop}\cite{O,GGP} \label{Cfti}
The Casson invariant is a degree $2$ finite type invariant.  
\end{prop}
So we can now use the theory of finite type invariants to study the behavior of the Casson invariant under Borromean surgery.  
\subsection{Proof of Theorem \ref{cy}} \label{proof2}
Let $L=L_1\cup ...\cup L_n$ be a disjoint union of Borromean surgery links in an integral homology sphere $M$.  
By using (\ref{MG}) and Proposition \ref{Cfti}, one can check that 
$$ \lambda(M_L) - \lambda(M) = \sum_{i=1}^{n}\left( \lambda(M_{L_i}) - \lambda(M) \right) + \sum_{1\le i<j\le n} \lambda([M;Y_i\cup Y_j]).  $$

So we are left with the computation of $\lambda([M;Y_i\cup Y_j])$ for all $ i<j$.  
By Theorem \ref{cy1}, we have 
\beQ
\lambda([M;Y_i\cup Y_j]) & = & \left( \lambda((M_{Y_j})_{Y_i}) - \lambda(M_{Y_j})\right) - 
\left( \lambda(M_{Y_i}) - \lambda(M)\right)\\
 & = & -2\mu_{123}(Y_i;M_{Y_j}) + 2\mu_{123}(Y_i;M),   
\eeQ
where we also use the fact that a Borromean surgery preserves the framings and linking numbers \cite{matv}.  
So, if we denote by $\Sigma_r$ a Seifert surface for $F^i_r$ ($r=1,2,3$), we see that $\lambda([M;Y_i\cup Y_j])$ counts, with signs, the triple intersection points between $\Sigma_1$, $\Sigma_2$ and $\Sigma_3$ created when doing surgery along $Y_j$.  

For each leaf $F^j_s$ of $Y_j$ ($s=1,2,3$), denote by $B_s$ the non-leaf component of $Y_j$ having linking number $1$ with $F^j_s$, oriented as in \S \ref{not}.  
Using Kirby moves, one sees that each positive (resp.\ negative) intersection point between $\Sigma_r$ and $F^j_s$ ($r,s=1,2,3$) contributes, \emph{via} surgery, to a band sum of $F^i_r$ with a parallel copy of $B_s$ with opposite (resp.\ same) orientation.  
So surgery along $Y_j$ results in band summing copies of $B_1$, $B_2$ and $B_3$ with the leaves of $Y_i$, with the number and orientation of these copies being completely determined by the linking numbers of the leaves of $Y_i$ and $Y_j$.  
Each nontrivial contribution to $\mu_{123}$ comes from a triple of band sums of pairwise distinct leaves of $Y_i$ with copies of pairwise distinct $B_s$'s.  
It thus follows from an easy counting argument that 
\[ \lambda([M;Y_i\cup Y_j]) = - 2 \sum_{\sigma\in  S_3} s(\sigma) l^{ij}_{1\sigma(1)}l^{ij}_{2\sigma(2)}l^{ij}_{3\sigma(3)},. \]  
which concludes the proof.  
\begin{rem} 
The theory of finite type invariants, that we used here to prove Theorem \ref{cy}, could actually also be used to give an alternative (and quite simple) proof of Theorem \ref{cy1}.  Such a proof, however, is only valid for integral homology spheres and is therefore not included in this paper.  
\end{rem}
\section{A remark and an application. } \label{rencom}
\subsection{Comparison to Lescop's surgery formula}
One can also use Lescop's global surgery formula to express the variation of the Casson invariant under a Borromean surgery,  
in terms of the multivariable Alexander polynomial $\Delta$ \cite{LG}.  

Consider the Borromean surgery link $L$ as in Figure \ref{conv} and denote by $B$ the standard Borromean link.  
Then \cite[1.4.8]{LG} gives that $\lambda(M_L) - \lambda(M)$ equals 
\[ -\left| \begin{array}{lll}
f_1 & l_{12} & l_{13} \\
l_{12} & f_{2} & l_{23} \\
l_{13} & l_{23} & f_{3} 
\end{array} \right|.\zeta(B) 
- \sum_{\circlearrowleft_{1,2,3}}\left| \begin{array}{ll}
f_1 & l_{12} \\
l_{12} & f_{2} 
\end{array} \right|.\zeta(B\cup F_3) 
- \sum_{\circlearrowleft_{1,2,3}} f_1.\zeta(B\cup F_2\cup F_3) - \zeta(L),\]
where, for an $n$-component link $K$, $\zeta(K)=\frac{\partial^n}{\partial t_1...\partial t_n}\Delta_K(1,...,1)$.  
For the Borromean link $B$ we know that $\zeta(B)=1$, and one can easily compute the values of $\zeta$ in
the two sums of the above formula: $\zeta(B\cup F_i)=0$ and $\zeta(B\cup F_i\cup F_j)=-l_{ij}$ for all $1\le i\ne j\le 3$.
It follows in particular that 
$$ \frac{\partial^6}{\partial t_1...\partial t_6}\Delta_L(1,...,1)=2.\mu_{123}(L) $$  
\subsection{An application.}
As a conclusion, we give an application of our main theorem.  Namely we show how 
Theorem \ref{cy1} provides an alternative proof of a realization theorem for the Casson invariant of homology spheres of Mazur type.  

Recall that a \emph{homology sphere of Mazur type} is obtained by surgery on $S^3$ along a $2$-component link $K_1\cup K_2$ 
such that $lk(K_1,K_2)=\pm 1$, $K_1$ is a trivial knot with framing $0$ and $K_2$ has framing $r\in \mathbf{Z}$.  
Such an integral homology sphere bounds a contractible $4$-manifold, and thus its Casson invariant is an even number.  
Mizuma showed the following.  
\begin{thm} \cite{mizuma}
For any even integer $k$, there is a homology sphere of Mazur type $M$ such that $\lambda(M)=k$.    
\end{thm}
Mizuma describes this homology sphere in terms of double branched cover of $S^3$ along so-called knots of 
\emph{$1$-fusion}.  
Theorem \ref{cy1} allows us to give explicit examples in terms of Borromean surgery links.  
Indeed, consider the link $L_{\pm n}$ depicted in the left part of Figure \ref{maz}.  
There, the three leaves have framing zero.  Clearly, by Theorem \ref{cy1}, surgery on $S^3$ along $L_{\pm n}$ produces an integral homology 
sphere with Casson invariant $\mp 2n$.  
\begin{figure}[!h] 
\begin{center} 
\input{mazurNEW.pstex_t} 
\caption{} \label{maz}
\end{center}
\end{figure} 
On the other hand, by Lemma \ref{KIRBY}, the link $L_{\pm n}$ is Kirby-equivalent to a $2$-component link $K_1\cup K_2$ as depicted on the 
right-hand side of the figure, which satisfies the above-mentioned conditions.  
Therefore $S^3_{L_{\pm n}}\cong S^3_{K_1\cup K_2}$ is a homology sphere of Mazur type.

\begin{thebibliography}{99} 
\newcommand{\no}{$\textrm{n}^{\circ}$} 
%
\bibitem{AL} E. Auclair, C. Lescop, \emph{Algebraic version of the calculus of clovers}, Algebr. Geom. Topol. \textbf{5} (2005), 
71-106.
%
\bibitem{GGP} S. Garoufalidis, M. Goussarov, M. Polyak, \emph{Calculus of clovers and finite type invariants of 3-manifolds}, Geom. Topol. $\mathbf{5}$ (2001), 75-108. 
%
\bibitem{G} M. Goussarov, \emph{Finite type invariants and $n$-equivalence of 3-manifolds}, Compt. Rend. Acad. Sc. Paris $\mathbf{329}$ S\'erie I (1999), 517-522. 
%
\bibitem{habi} K. Habiro, \emph{Claspers and finite type invariants of links}, Geom. Topol. $\mathbf{4}$ (2000), 1--83. 
%
\bibitem{johannes} J. Johannes, \emph{A Type 2 Polynomial Invariant of Links Derived from the Casson-Walker Invariant},  J. Knot Theory Ram.  \textbf{8} (1999) 491-504. 
%
\bibitem{LG} C. Lescop, \emph{Global surgery formula for the Casson-Walker invariant}, Annals of Math. Studies \textbf{140}, Princeton Univ. Press (1996).
%
\bibitem{L} C. Lescop, \emph{A sum formula for the Casson-Walker invariant}, Invent. Math., \textbf{133} (1998) 613-681. 		
%
\bibitem{massuyeau} G. Massuyeau, \emph{Spin Borromean surgeries}, Trans. Amer. Math. Soc. $\mathbf{355}$ (2003), 3991-4017. 
%
\bibitem{matv} S. Matveev, \emph{Generalized surgery of three-dimensional manifolds	and representation of homology spheres}, Math. Notices Acad. Sci. \textbf{42:2} (1988), 651-656.  
%
\bibitem{mizuma} Y. Mizuma, \emph{On the Casson invariant of Mazur's homology spheres}, preprint.  
%
\bibitem{O} T. Ohtsuki, \emph{Finite type invariants of integral homology $3$-spheres}, J. Knot Theory Ram. \textbf{5} (1996), 101-115.
%
\bibitem{W} K. Walker,  \emph{An Extension of Casson's Invariant}, Annals of Math. Studies \textbf{126}, Princeton Univ. Press (1992). 
%
\end{thebibliography}
\end{document}